\documentclass[11pt]{article}

\usepackage{amssymb,amsmath,bm}
\usepackage{amsthm}

\usepackage{textcomp}
\usepackage{enumerate}      
\usepackage{graphicx}        
\usepackage{caption}

\usepackage{mathrsfs}

\usepackage{url}

\renewcommand{\setminus}{{\smallsetminus}}

\newtheorem*{namedtheorem}{\theoremname}
\newcommand{\theoremname}{testing}

\newtheorem{theorem}{Theorem}[section]

\newtheorem{definition}[theorem]{Definition}

\newtheorem{conjecture}[theorem]{Conjecture}

\theoremstyle{remark}

\theoremstyle{remark}

\numberwithin{equation}{section}

\makeatother
\DeclareMathOperator{\Z}{\mathbb{Z}}

\DeclareMathOperator{\vol}{\mathrm{Vol}}

\begin{document}
\title{Growth of Turaev-Viro invariants and cabling} 
\author{Renaud Detcherry} \date{}
\maketitle

\begin{abstract} 
The Chen-Yang volume conjecture \cite{Chen-Yang} states that the growth rate of the Turaev-Viro invariants of a compact oriented $3$-manifold determines its simplicial volume. In this paper we prove that the Chen-Yang conjecture is stable under $(2n+1,2)$ cabling.
\end{abstract}

\section{Introduction}
For $M$ a $3$-manifold that is either closed or with boundary, the Turaev-Viro invariants $TV_r(M)$ are real valued topological invariants that can be computed using state sums over a triangulation of $M.$ They depend of the choice of an $2r$-th root of unity $q$ whose square is a primitive $r$-th root; in this paper we will always choose $q=e^{\frac{2i\pi}{r}}$ and assume $r$ odd. Moreover, when $M$ is a manifold with empty or toroidal boundary, the invariants $TV_r(M)$ are always non-negative.
\\ The geometric meaning of the Turaev-Viro invariants is hard to understand from their state sum definition. However, a conjecture of Chen and Yang \cite{Chen-Yang} states that the asymptotics of the $TV_r$ invariants at root $q=e^{\frac{2i\pi}{r}}$ is related to hyperbolic volume: 
\begin{conjecture}\label{TVVC}\cite{Chen-Yang}
For any hyperbolic manifold $M$ (closed or with boundary), we have
$$\underset{r \rightarrow \infty, \ r \ \textrm{odd}}{\lim} \frac{2\pi}{r}\log |TV_r(M,q=e^{\frac{2i\pi}{r}})|=\vol (M)$$
where $\vol (M)$ is the hyperbolic volume of $M.$
\end{conjecture}
 Conjecture \ref{TVVC} is reminiscent of the Volume Conjecture of Kashaev and Murakami-Murakami \cite{kashaev:vol-conj}\cite{murakami:vol-conj} where the $TV_r$ invariants are replaced with evaluations $J_n(K,e^{\frac{2i\pi}{n}})$ of the normalized colored Jones polynomials of an hyperbolic knot.  
\\ In \cite{DKY}, Yang, Kalfagianni and the author gave a formula relating the Turaev-Viro invariants of a link complement to colored Jones polynomials, establishing a connection between the two conjectures. In the same paper, Conjecture \ref{TVVC} was proved for the complements of the figure-eight knot, the borromean link and knots of Gromov norm zero. 
\\
\\ Let us define the growth rate of Turaev-Viro invariants by:
\begin{definition}Let $M$ be a $3$-manifold, closed or with boundary. Then the Turaev-Viro growth rate is
$$LTV(M)=\underset{r \rightarrow \infty, \ r \ \textrm{odd}}{\limsup} \frac{2\pi}{r}\log |TV_r(M,q=e^{\frac{2i\pi}{r}})|.$$
\end{definition}
A way to restate Conjecture \ref{TVVC} for general $3$-manifold (not necessarily hyperbolic), is the following:
\begin{conjecture}\label{TVVCgeneral} For any compact oriented $3$-manifold $M,$ we have
$$LTV(M)=\vol(M),$$
where $\vol(M)$ is the simplicial volume of $M.$
\end{conjecture}
We recall that the simplicial volume of $M$ can be thought either as the sum of the hyperbolic volumes of the hyperbolic pieces in the JSJ decomposition of $M,$ or as $v_3||M||$ where $||M||$ is the Gromov norm of $M.$ Note that the Turaev-Viro invariants sometimes vanish on lens spaces, thus replacing the limit by an upper limit is necessary.
\\ In \cite{DK}, Kalfagianni and the author investigated the growth rate of Turaev-Viro invariants. They showed that the growth rate of Turaev-Viro invariants satisfies properties similar to that of the simplicial volume:
\begin{theorem}\label{thm:LTVprop}\cite{DK} Let $M$ be a compact oriented $3$-manifold, with empty or toroidal boundary.
\begin{itemize}
\item[(1)] If $M$ is a Seifert manifold, then there exists constants $B>0$ and $N$ such that for any odd $r\geqslant 3,$ we have $TV_r(M)\leqslant B r^N$ and $LTV(M)\leqslant 0.$
\item[(2)] If $M$ is a Dehn-filling of $M',$ then $TV_r(M)\leqslant TV_r(M')$ and $LTV(M)\leqslant LTV(M').$
\item[(3)] If $M=M_1 \underset{T}{\cup} M_2$ is obtained by gluing two $3$-manifolds $M_1$ and $M_2$ along a torus boundary component, then $TV_r(M)\leqslant TV_r(M_1)TV_r(M_2)$ and $LTV(M)\leqslant LTV(M_1)+LTV(M_2).$
\end{itemize}
\end{theorem}
These properties are parallel to the properties of the simplicial volume: the simplicial volume of Seifert manifolds is $0,$ the simplicial volume decreases under Dehn-filling and is subadditive under gluing along tori \cite{thurston:notes}.
\\
\\ Let $p,q$ be coprime integers, the $(p,q)$-cabling space is the complement of a $(p,q)$-torus knot standardly embedded in a solid torus. A $(p,q)$-cabling of a manifold $M$ with toroidal boundary is a manifold $M'$ obtained by gluing a $(p,q)$-cabling space to a boundary component of $M.$  In this paper, we will investigate the compatibility of Conjecture \ref{TVVC} with $(p,2)$-cabling. We will show the following:
\begin{theorem}\label{thm:(p,2)cabling}
Let $M$ be a manifold with toroidal boundary and $M'$ be obtained by gluing a $(p,2)$-cabling space $C_{p,2}$ to a boundary component of $M.$ Then there are constants $B>0$ and $N>0$ such that $$\frac{1}{Br^N}TV_r(M')\leqslant TV_r(M)\leqslant Br^N TV_r(M').$$
\end{theorem}
In particular, this means that if Conjecture \ref{TVVC} is true for $M$ then it is true for $M'.$
\section{Preliminaries}
\label{sec:prelim}

\subsection{Reshetikhin-Turaev $\mathrm{SO}_3$-TQFTs and TQFT basis}
We briefly sketch the properties of Reshetikhin-Turaev $\mathrm{SO}_3$-TQFTs, defined by Reshetikhin and Turaev in \cite{ReTu}. We will introduce them in the skein-theoretic framework of Blanchet, Habegger, Masbaum and Vogel \cite{BHMV2}. We refer to \cite{BHMV1}\cite{BHMV2} for the details of these constructions.
\\
\\ For any odd integer $r\geqslant 3,$ and primitive $2r$-th root of unity $A,$ there is an associated TQFT functor $RT_r,$ with the following properties:
\begin{itemize}
\item[-]For $\Sigma$ a closed compact oriented surface,  $RT_r(\Sigma)$ is a finite dimensional $\mathbb{C}$-vector space, with a natural Hermitian form. Moreover for a disjoint union $\Sigma \coprod \Sigma'$ one has $RT_r(\Sigma \coprod \Sigma')=RT_r(\Sigma) \otimes RT_r(\Sigma').$
\item[-]For $M$ a compact oriented closed $3$-manifold, $RT_r(M)$ is the $\mathrm{SO}_3$ Reshetikhin-Turaev invariant, a complex valued topological invariant of $3$-manifolds, and for $M$ with boundary, $RT_r(M)$ is a vector in $RT_r(\partial M).$
\item[-]If $(M,\Sigma_1,\Sigma_2)$ is a cobordism, $RT_r(M):RT_r(\Sigma_1) \rightarrow RT_r(\Sigma_2)$ is a linear map. Moreover, the composition of cobordisms is sent to the composition of linear maps, up to a power of $A.$
\end{itemize}
Moreover, some basis of the TQFT spaces $RT_r(\Sigma)$ of surfaces has been explicitely described in \cite{BHMV2}. We recall that for $\Sigma$ a surface, $RT_r(\Sigma)$ is a quotient of the Kauffman module of a handlebody of the same genus. 
\\ In the case of a torus we get the following picture: the torus $T^2$ is the boundary of a solid torus $D^2\times S^1.$ One gets a family of elements of the Kauffman module of $D^2\times S^1$ by taking the core $\lbrace 0 \rbrace \times S^1$ and coloring it by the $i-1$-th Jones-Wenzl idempotents, thus obtaining an element $e_i \in RT_r(T^2).$ For a definition of Jones-Wenzl idempotents we refer to \cite{BHMV2}. For $r=2m+1,$ and $A$ a $2r$-th root of unity, only finitely many Jones-Wenzl idemptotents can be defined, thus only the elements $e_1,\ldots ,e_{2m-1}$ are well defined (see \cite{BHMV2}[Lemma 3.2]).
\\ As elements of the Kauffman module of the solid torus, the $e_i$'s can also be considered as elements of $RT_r(T^2).$ One gets a basis of $RT_r(T^2)$ consisting of elements $e_i:$
\begin{theorem}\label{thm:basis}\cite{BHMV2}[Theorem 4.10]
If $r=2m+1\geqslant 3,$ then the family $e_1,e_2,\ldots e_{m}$ is an orthonormal basis of $RT_r(T^2).$ Moreover one has $e_{m-i}=e_{m+1+i}$ for $0\leqslant i \leqslant m-1.$
\end{theorem}
Note that the last part implies that the family $e_1,e_3,\ldots e_{2m-1}$ is the same basis of $RT_r(T^2)$ as the family $e_1,e_2,\ldots e_m,$ in a different order.
\\ 
\\ As a consequence of TQFTs axioms, the $RT_r$ vectors associated to link complements can be tied to values of colored Jones polynomials. Indeed, if $M=S^3\setminus L$ is a link complement where $L$ has $n$ components, $RT_r(M)$ will be a vector in $RT_r(T^2)^{\otimes n}$ whose coefficient along $e_{i_1}\otimes e_{i_2} \ldots \otimes e_{i_n}$ is obtained by gluing solid tori with cores colored by Jones-Wenzl idempotents to the boundary components of $M$ and taking the Reshetikhin-Turaev invariant of that. Thus the coefficient we get is $\eta_r J_{i}(L,A^4),$ where $J_{i}(L,t)$ is the $i$-th unnormalized colored Jones polynomial of $L,$ $i$ is a multi-index of colors, and 
$$\eta_r=RT_r(S^3)=\frac{A^2-A^{-2}}{\sqrt{-r}}.$$
\subsection{Relationship with the Turaev-Viro invariants}
While the Turaev-Viro invariants of compact oriented $3$-manifolds $M$ are defined as state sums over triangulations of $M$ (see \cite{TuraevViro}), we will only use a well-known identity relating the Turaev-Viro invariants and Reshetikhin-Turaev invariants. This property was first proved by Roberts \cite{Roberts} in the case of closed $3$-manifolds, and extended to manifolds with boundary by Benedetti and Petronio  \cite{BePe}. For simplicity we state it only in the case of manifolds with toroidal boundary:
\begin{theorem}\label{thm:TV-RT}
Let $M$ be a compact oriented manifold with toroidal boundary, let $r\geqslant 3$ be an odd integer and let $A$ be a primitive $2r$-th root of unity. Then have 
$$TV_r(M,A^2)=||RT_r(M,A)||^2$$
where $|| \cdot ||$ is the natural Hermitian norm on $RT_r(\partial M).$ 
\end{theorem}
For $T^2$ a torus, the natural Hermitian form on $RT_r(T^2)$ is definite positive for any $A,$ thus $TV_r(M)$ is non-negative.
\subsection{The cabling formula}
For $r=2m+1,$ it is convenient to extend the definition of vectors $e_i \in RT_r(T^2)$ to all possible values of $i \in \Z$ as follows: we formally set $e_{-l}=-e_{l}$ for any $l \geqslant 0$ (and particular $e_0=0$) and 
$$e_{l+kr}=(-1)^k e_{l}$$ for any $k \in \mathbb{Z}.$ This is compatible with the above mentioned symmetry of the vectors $e_i$'s.
\\ 
\\ Recall that for $p$ and $q$ coprime integers with $q>0,$ the $(p,q)$-cabling space $C_{p,q}$ is the complement in a solid torus of a standardly embedded $(p,q)$-torus knot. The cobordism $C_{p,q}$ gives rise in TQFT to a linear map
$$RT_r(C_{p,q}) : RT_r(T^2) \rightarrow RT_r(T^2).$$
We now describe the action of this map on the basis $\lbrace e_i \rbrace_{1 \leqslant i \leqslant m}$ of $RT_r(T^2).$ By TQFT axioms, the map $RT_r(C_{p,q})$ sends the element $e_i$ to the element of $RT_r(T^2)$ corresponding to a $(p,q)$-torus knot embedded in the solid torus and colored by the $i-1$-th Jones-Wenzl idempotent.
\\ Morton computed these elements using skein calculus, yielding the following formula:
\begin{theorem}\cite{Morton} \label{thm:cablingformula}(Cabling formula)
\\ For any odd $r=2m+1 \geqslant 3,$ for any $1 \leqslant i \leqslant m,$ for any coprime integers $(p,q),$ one has:
$$RT_r(C_{p,2})(e_i)=A^{pq(i^2-1)/2}\underset{k\in \mathcal{S}_i}{\sum} A^{-2pk(qk+1)} e_{2qk+1}.$$
where $\mathcal{S}_i$ is the set
$$\mathcal{S}_i =\lbrace -\frac{i-1}{2}, -\frac{i-3}{2},\ldots \frac{i-3}{2},\frac{i-1}{2}\rbrace.$$ 
\end{theorem} 
\section{Stability under $(p,2)$-cabling}
\label{sec:(p,2)cabling}
In this section, we will let $r=2m+1 \geqslant 3$ be an odd integer and let $A=e^{\frac{i\pi}{r}},$ which is a primitive $2r$-th root of unity.
\\
\\
Recall that a $(p,q)$-cabling $M'$ of a manifold $M$ consist of gluing the exterior torus boundary component of $C_{p,q}$ to a boundary component of $M.$ The JSJ decomposition of $M'$ will consist of the pieces in the JSJ decomposition of $M,$ plus an extra piece that is the cabling space $C_{p,q}.$ As $C_{p,q}$ is a Seifert manifold, $M$ and $M'$ have the same simplicial volume.
\begin{theorem}\label{thm:(p,2)cabling}
Let $M$ be a manifold with toroidal boundary, let $p$ be an odd integer and let $M'$ be a $(p,2)$-cabling of $M.$ Then there exists constants $B>0$ and $N$ such that
$$\frac{r^{-N}}{B} TV_r(M) \leqslant TV_r(M') \leqslant B r^N TV_r(M)$$
In particular, we have $LTV(M)=LTV(M').$
\end{theorem}
As $M$ and $M'$ in the theorem have the same volume, the conclusion implies that if Conjecture \ref{TVVC} is true for $M$ it is true for $M'.$
\begin{proof}
First, $M'$ is obtained by gluing $C_{p,2}$, which is Seifert and thus has volume $0,$ to $M$ along a torus. We know by Theorem \ref{thm:LTVprop} that 
$$TV_r(M')\leqslant TV_r(M)TV_r(C_{p,2}).$$
But as $C_{p,2}$ has volume $0,$ one has $TV_r(C_{p,2})\leqslant B r^N$ for some constants $B>0$ and $N.$
\\ To prove the other inequality, we study the map induced by $C_{p,2}$ in the $RT_r$-TQFT. If $T$ is the boundary component coming from the $(p,2)$-torus knot and $T'$ is the exterior boundary component, then $C_{p,2}$ induces a map
$$RT_r(C_{p,2}):RT_r(T)\rightarrow RT_r(T').$$
If $M$ has only one boundary component, as $RT_r$ is a TQFT, we have that $RT_r(M')=RT_r(C_{p,2})(RT_r(M)).$ If $M$ has other boundary components than the one used to glue $C_{p,2},$ then for any coloring $i$ of the other boundary components of $M$, we have
$$RT_r(M',i)=RT_r(C_{p,2})(RT_r(M,i)).$$
In all cases, if $RT_r(C_{p,2})$ is invertible, we can write 
$$||RT_r(M)||\leqslant |||RT_r(C_{p,2})^{-1}||| \cdot ||RT_r(M')||,$$
where $|| \cdot ||$ is the norm induced by the Hermitian form on $RT_r(\partial M)$ or $RT_r(\partial M' )$ and $||| \cdot |||$ is the corresponding operator norm.
\\ To conclude the proof of the theorem, we thus only need to prove that $RT_r(C_{p,2})$ is invertible, and that $|||RT_r(C_{p,2})^{-1}|||$ grows at most polynomially.
\\ We can compute the matrix of $RT_r(C_{p,2})$ in the basis $e_1,\ldots e_{m}$ of $RT_r(T)$ by the cabling formula recalled above as Theorem \ref{thm:cablingformula}.
\\ For $q=2,$ the formula states:
$$RT_r(C_{p,2})(e_i)=A^{p(i^2-1)}\underset{k \in \mathcal{S}_i}{\sum} A^{-2pk(2k+1)} e_{4k+1}$$
The cabling formula implies that the image lies in the vector space spanned by $e_1,e_3,\ldots , e_{2m-1},$ but by the Symmetry Principle recalled in \ref{thm:basis}, $e_{m-i}=e_{m+1+i}$ for all $0 \leqslant i \leqslant m-1.$ Thus $\lbrace e_1,e_3,\ldots , e_{2m-1} \rbrace$ is actually the basis $\lbrace e_1,\ldots e_{m} \rbrace$ in a different order.
\\ From the cabling formula we get that $RT_r(C_{p,2})(e_i)$ lies in $\mathrm{Span}(e_1,e_3,\ldots e_{2i-1})$ and that the coefficient in $e_{2j-1}$ in $RT_r(C_{p,2})(e_i)$ is $A^{p(i^2-1)} \cdot A^{p(j-j^2)}$ if $j$ has same parity as $i$, and $-A^{p(i^2-1)} \cdot A^{p(j-j^2)}$ else. We can write $RT_r(C_{p,2})$ in the basis $\lbrace e_1,\ldots e_{m} \rbrace$ and $\lbrace e_1,e_3,\ldots , e_{2m-1} \rbrace$ as a product of two diagonal matrices and a triangular matrix:
$$RT_r(C_{p,2})=\begin{pmatrix}
1 & 0 & \ldots & 0
\\ 0 & A^{p(2-2^2)} & \ddots  & \vdots
\\  \vdots & \ddots & \ddots & 0
\\ 0 & \ldots  & 0 & A^{p(m-m^2)}
\end{pmatrix}
\begin{pmatrix}1 & -1 & 1 &\ldots & 
\\ 0 & 1 & -1 & & 
\\  \vdots & &  1 & \ddots & \vdots
\\ & & & \ddots & -1
\\ 0 &  & \ldots & &  1

\end{pmatrix}
\begin{pmatrix}
1 & 0 & \ldots & 0
\\ 0 & A^{p(2^2-1)} & \ddots  & \vdots
\\  \vdots & \ddots & \ddots & 0
\\ 0 & \ldots  & 0 & A^{p(m^2-1)}
\end{pmatrix}
$$
And thus we have:
$$RT_r(C_{p,2})^{-1}=\begin{pmatrix}
1 & 0 & \ldots & 0
\\ 0 & A^{p(1-2^2)} & \ddots  & \vdots
\\  \vdots & \ddots & \ddots & 0
\\ 0 & \ldots  & 0 & A^{p(1-m^2)}
\end{pmatrix}
\begin{pmatrix}1 & 1 & 0 &\ldots & 0 
\\ 0 & 1 & 1 & & \vdots
\\  \vdots & \ddots &  1 & \ddots & 
\\ & & & \ddots & 1
\\ 0 &  & \ldots & 0 &  1

\end{pmatrix}
\begin{pmatrix}
1 & 0 & \ldots & 0
\\ 0 & A^{p(2^2-2)} & \ddots  & \vdots
\\  \vdots & \ddots & \ddots & 0
\\ 0 & \ldots  & 0 & A^{p(m^2-m)}
\end{pmatrix}
$$
The two diagonal matrices are isometries and the middle matrix clearly has norm bounded by a polynomial in $r,$ which concludes the proof of Theorem \ref{thm:(p,2)cabling}.
\end{proof}
\bibliographystyle{hamsplain}
\bibliography{biblio}

 \noindent 
Renaud Detcherry\\
Department of Mathematics, Michigan State University\\
East Lansing, MI 48824\\
(detcherry@math.msu.edu)

\end{document}